\documentclass[12pt]{article}
 \usepackage{latexsym}
 \usepackage{amssymb}
 \usepackage{graphicx}

 \newtheorem{Theorem}{Theorem}

\newtheorem{Proposition}{Proposition}

\newtheorem{Algorithm}{Algorithm}
\newtheorem{Question}{Question}

\newcommand{\A}{{\cal A}}
\newcommand{\B}{{\cal B}}

\newcommand{\PP}{{\cal P}}

\newcommand{\uu}{{\bf u}}
\newcommand{\vv}{{\bf v}}
\newcommand{\x}{{\bf x}}
\newcommand{\y}{{\bf y}}
\newcommand{\z}{{\bf z}}
\newcommand{\e}{{\bf e}}
\newcommand{\0}{{\bf 0}}
\newcommand{\w}{{\bf w}}

\newcommand{\qed}{\nobreak \ifvmode \relax \else
      \ifdim\lastskip<1.5em \hskip-\lastskip
      \hskip1.5em plus0em minus0.5em \fi \nobreak
      \vrule height0.75em width0.5em depth0.25em\fi}

\def \ep{\hbox{ }\hfill$\Box$}

\begin{document}
\title{Hankel Tensors: Associated Hankel Matrices and Vandermonde Decomposition}

\author{
Liqun Qi \thanks{Email: maqilq@polyu.edu.hk. Department of Applied
Mathematics, The Hong Kong Polytechnic University, Hung Hom,
Kowloon, Hong Kong. This author's work was supported by the Hong
Kong Research Grant Council (Grant No. PolyU 501909, 502510, 502111
and 501212).}}

\date{\today} \maketitle

\begin{abstract}
\noindent  
Hankel tensors arise from applications such as signal processing. In
this paper, we make an initial study on Hankel tensors.   For each
Hankel tensor, we associate it with a Hankel matrix and a higher
order two-dimensional symmetric tensor, which we call the associated
plane tensor.   If the associated Hankel matrix is positive
semi-definite, we call such a Hankel tensor a strong Hankel tensor.
We show that an $m$ order $n$-dimensional tensor is a Hankel tensor
if and only if it has a Vandermonde decomposition.   We call a
Hankel tensor a complete Hankel tensor if it has a Vandermonde
decomposition with positive coefficients.   We prove that if a
Hankel tensor is copositive or an even order Hankel tensor is
positive semi-definite, then the associated plane tensor is
copositive or positive semi-definite, respectively.   We show that
even order strong and complete Hankel tensors are positive
semi-definite, the Hadamard product of two strong Hankel tensors is
a strong Hankel tensor, and the Hadamard product of two complete
Hankel tensors is a complete Hankel tensor. We show that all the
H-eigenvalue of a complete Hankel tensors (maybe of odd order) are
nonnegative. We give some upper bounds and lower bounds for the
smallest and the largest Z-eigenvalues of a Hankel tensor,
respectively.   Further questions on Hankel tensors are raised.

\noindent {\bf Key words:}\hspace{2mm} Hankel tensors, Hankel
matrices, plane tensors, positive semi-definiteness,
co-positiveness, generating functions, Vandermonde decomposition,
eigenvalues of tensors. \vspace{3mm}

\noindent {\bf AMS subject classifications (2010):}\hspace{2mm}
15A18; 15A69
  \vspace{3mm}

\end{abstract}


\section{Introduction}
\hspace{4mm} Hankel matrices play an important role in linear
algebra and its applications \cite{Ge, HJ, Wid}.  As a natural
extension of Hankel matrices, Hankel tensors arise from applications
such as signal processing.

Denote $[n] := \{ 1, \cdots, n \}$.  Let $\A = \left(a_{i_1 \cdots
i_m}\right)$ be a real $m$th order $n$-dimensional tensor. If there
is a vector $\vv = (v_0, v_1, \cdots, v_{(n-1)m})^\top$ such that
for $i_1, \cdots, i_m \in [n]$, we have
\begin{equation}\label{e1}
a_{i_1\cdots i_m} \equiv v_{i_1+i_2+\cdots +i_m-m}, \end{equation}
then we say that $\A$ is an $m$th order {\bf Hankel tensor}. Hankel
tensors were introduced by Papy, De Lathauwer and Van Huffel in
\cite{PDV} in the context of the harmonic retrieval problem, which
is at the heart of many signal processing problems.  In \cite{BB},
Badeau and Boyer proposed fast higher-order singular value
decomposition (HOSVD) for third order Hankel tensors.

A real $m$th order $n$-dimensional tensor (hypermatrix) $\A =
(a_{i_1\cdots i_m})$ is a multi-array of real entries $a_{i_1\cdots
i_m}$, where $i_j \in [n]$ for $j \in [m]$. Denote the set of all
real $m$th order $n$-dimensional tensors by $T_{m, n}$.  Then $T_{m,
n}$ is a linear space of dimension $n^m$.    If the entries
$a_{i_1\cdots i_m}$ are invariant under any permutation of their
indices, then $\A$ is a {\bf symmetric tensor}. Denote the set of
all real $m$th order $n$-dimensional symmetric tensors by $S_{m,
n}$. Then $S_{m, n}$ is a linear subspace of $T_{m, n}$. Clearly, a
Hankel tensor is a symmetric tensor. Denote the set of all real
$m$th order $n$-dimensional Hankel tensors by $H_{m, n}$. Then
$H_{m, n}$ is a linear subspace of $S_{m, n}$, with dimension
$(n-1)m+1$.

Let $\A = (a_{i_1\cdots i_m}) \in S_{m, n}$ and $\x = (x_1, \cdots,
x_n)^\top \in \Re^n$. Denote $$\A \x^m = \sum_{i_1,\cdots, i_m=1}^n
a_{i_1\cdots i_m}x_{i_1}\cdots x_{i_m}.$$ Denote $\Re^n_+ = \{ \x
\in \Re^n : \x \ge \0 \}$.   If $\A \x^m \ge 0$ for all $\x \in
\Re^n_+$, then $\A$ is called {\bf copositive}.  If $\A \x^m > 0$
for all $\x \in \Re^n_+, \x \not = \0$, then $\A$ is called {\bf
strongly copositive} \cite{Qi0}. Suppose that $m$ is even. If $\A
\x^m \ge 0$ for all $\x \in \Re^n$, then $\A$ is called {\bf
positive semi-definite}.  If $\A \x^m > 0$ for all $\x \in \Re^n, \x
\not = \0$, then $\A$ is called {\bf positive definite} \cite{Qi}.
Positive semi-definite symmetric tensors are useful in automatical
control \cite{Qi} and higher-order diffusion tensor imaging
\cite{BV, CDHS, HHNQ, QYW, QYX}.  It is established in \cite{Qi}
that an even order symmetric tensor $\A \in S_{m, n}$ is positive
semi-definite if and only if all of its H-eigenvalues (or
Z-eigenvalues) are nonnegative.   On the other hand, copositive
tensors do not restrict the order to be even, thus are more general.
Nonnegative tensors, positive semi-definite tensors and Laplacian
tensors \cite{Qi1} are copositive tensors \cite{Qi0}.

In the next section, for each Hankel tensor $\A \in H_{m, n}$, we
associate it with a symmetric tensor $\PP \in S_{(n-1)m, 2}$. We
call such a tensor the {\bf associated plane tensor}.   We use the
term ``plane tensor'' here as its dimension is only $2$,
corresponding to a tensor on the plane in physics, while three
dimensional tensors are called space tensors in \cite{QY}. Actually,
$S_{l, 2} \equiv H_{l, 2}$ for any $l \ge 2$.  But we do not stress
that $\PP$ is a Hankel tensor here.  For a symmetric tensor, we may
use the elimination method proposed in \cite{QWW} to calculate its
Z-eigenvalues, and to determine if it is positive semi-definite or
not when the order is even.  We show that if a Hankel tensor is
copositive or an even order Hankel tensor is positive semi-definite,
then the associated plane tensor is copositive or positive
semi-definite, respectively.

Suppose that $\A \in H_{m, n}$ is defined by (\ref{e1}).    Let $A =
(a_{ij})$ be an $\lceil {(n-1)m+2 \over 2} \rceil  \times \lceil
{(n-1)m+2 \over 2} \rceil$ matrix with $a_{ij} \equiv v_{i+j-2}$,
where $v_{2\lceil {(n-1)m \over 2} \rceil}$ is an additional number
when $(n-1)m$ is odd. Then $A$ is a Hankel matrix, associated with
the Hankel tensor $\A$. Such an associated Hankel matrix is unique
if $(n-1)m$ is even. If the Hankel matrix $A$ is positive
semi-definite, then we say that $\A$ is a {\bf strong Hankel
tensor}.

It is clear that the Hadamard product of two Hankel tensors is a
Hankel tensor.   In Section 3, we show that an even order strong
Hankel tensor is positive semi-definite and the Hadamard product of
two strong Hankel tensors is also a strong Hankel tensor. In order
to do this, we introduce a generating function for a Hankel tensor.
We show that a Hankel tensor has a nonnegative generating function
if and only if it is a strong Hankel tensor.  We give an example of
a positive semi-definite Hankel tensor, which is not a strong Hankel
tensor, and an example that the Hadamard product of two positive
semi-definite Hankel tensors is not positive semi-definite.

In Section 4, we introduce Vandermonde decomposition and show that
an $m$ order $n$-dimensional tensor is a Hankel tensor if and only
if it has a Vandermonde decomposition.   We call a Hankel tensor a
{\bf complete Hankel tensor} if it has a Vandermonde decomposition
with positive coefficients.   We show that an even order complete
Hankel tensor is positive semi-definite and the Hadamard product of
two complete Hankel tensors is also a complete Hankel tensor..   In
general, a positive semi-definite Hankel tensor may not be a
complete Hankel tensor.

As even order complete and strong Hankel tensors are positive
semi-definite symmetric tensors, all of their H-eigenvalues and
Z-eigenvalues are nonnegative, by Theorem 5 of \cite{Qi}.  On the
other hand, what are the spectral properties of odd order complete
and strong Hankel tensors?  We study these in Section 5. We show
that all the H-eigenvalue of an odd order complete Hankel tensors
are also nonnegative.   Suppose that $\x = (x_1, \cdots, x_n)^\top$
is a Z-eigenvector of a complete or strong Hankel tensor $\A$,
associated with a nonzero Z-eigenvalue $\lambda$ of $\A$.   We show
that for all odd $i$, $x_i \ge 0$ if $\lambda > 0$, and $x_i \le 0$
if $\lambda < 0$.   If $\A$ is a complete Hankel tensor, then $x_1 >
0$ if $\lambda > 0$, and $x_1 < 0$ if $\lambda < 0$.

In Section 6, we give some upper bounds and lower bounds for the
smallest and the largest Z-eigenvalues of a Hankel tensor,
respectively.   In Section 7, we present an algorithm to determine a
symmetric plane tensor $\PP \in S_{l, 2}$ is copositive or not for
$l \ge 2$.

Several questions are raised in Sections 2-6.  Some further
questions are raised in Section 8.

Throughout this paper, we assume that $m, n \ge 2$.   We use small
letters $x, u, v, \alpha, \cdots$, for scalers, small bold letters
$\x, \y, \uu, \cdots$, for vectors, capital letters $A, B, \cdots$,
for matrices, calligraphic letters $\A, \B, \cdots$, for tensors.
Denote $\e_i \in \Re^n$ as the $i$th unit vector for $i \in [n]$,
and $\0$ as the zero vector in $\Re^n$.

\section{Associated Plane Tensors, Copositive Hankel Tensors, Positive Semi-Definite Hankel Tensors}

We first give a necessary condition for a Hankel tensor to be
copositive.

\begin{Proposition}  \label{p1}
Suppose that $\A \in H_{m, n}$ is defined by (\ref{e1}). If $\A$ is
copositive, then $v_{(i-1)m} \ge 0$ for $i \in [n]$.
\end{Proposition}
\noindent {\bf Proof.} Since $v_{(i-1)m} = \A (\e_i)^m$ for $i \in
[n]$, the conclusion follows from the definition of copositive
tensors. \ep

As a positive semi-definite symmetric tensor is copositive
\cite{Qi0}, the condition $v_{(i-1)m} \ge 0$ for $i \in [n]$ is also
a necessary condition for an even order Hankel tensor to be positive
semi-definite.

\medskip

For any nonnegative integer $k$, define $s(k, m, n)$ as the number
of distinct ordered sets of indices $(i_1,\cdots, i_m)$ such that
$i_j \in [n]$ for $j \in [m]$ and $i_1+ \cdots + i_m - m = k$. Then
$s(0, m, n) = 1, s(1, m, n) = m, s(2, m, n) = {m(m+1) \over 2},
\cdots$.

We now define the associated plane tensor of a Hankel tensor.
Suppose that $\A \in H_{m, n}$ is defined by (\ref{e1}).  Define
$\PP = (p_{i_1\cdots i_{(n-1)m}}) \in S_{(n-1)m, 2}$ by
$$p_{i_1\cdots i_{(n-1)m}} = { s(k, m, n) v_k \over \left({(n-1)m \atop
k}\right)},$$ where $k = i_1+  \cdots + i_{(n-1)m} - (n-1)m$.    We
call $\PP$ the {\bf associated plane tensor} of $\A$.

\begin{Theorem} \label{t0}
If a Hankel tensor $\A \in H_{m, n}$ is copositive, then its
associated plane tensor $\PP$ is copositive.  If an even order
Hankel tensor $\A \in H_{m, n}$ is positive semi-definite, then its
associated plane tensor $\PP$ is positive semi-definite.
\end{Theorem}
\noindent {\bf Proof.} Suppose that $\A$ is copositive.  By
Proposition \ref{p1}, $v_{(n-1)m} \ge 0$.   Let $\y = (y_1,
y_2)^\top \in \Re^2_+$.  If $y_1 = y_2 = 0$, then clearly $\PP
\y^{(n-1)m} = 0$. If $y_1 = 0$ and $y_2 \not = 0$, then $\PP
\y^{(n-1)m} =v_{(n-1)m}y_2^{(n-1)m} \ge 0$.  We now assume that $y_1
\not = 0$. Let $u = {y_2 \over y_1}$.  Then $u \ge 0$.   We have
\begin{equation} \label{e0}
\PP \y^{(n-1)m} = y_1^{(n-1)m} \sum_{k=0}^{(n-1)m} \left({(n-1)m
\atop k}\right) \cdot { s(k, m, n) v_k \over \left({(n-1)m \atop
k}\right)} u^k = y_1^{(n-1)m} \A \uu^m \ge 0,
\end{equation} where $\uu = (1, u, u^2, \cdots, u^{n-1})^\top \in
\Re^n_+$.  Thus, $\PP$ is copositive.

Suppose that $m$ is even and $\A$ is positive semi-definite.   Then
$(n-1)m$ is also even.  By Proposition \ref{p1}, $v_{(n-1)m} \ge 0$.
Let $\y = (y_1, y_2)^\top \in \Re^2$.  If $y_1 = y_2 = 0$, then
clearly $\PP \y^{(n-1)m} = 0$. If $y_1 = 0$ and $y_2 \not = 0$, then
$\PP \y^{(n-1)m} =v_{(n-1)m}y_2^{(n-1)m} \ge 0$.  We now assume that
$y_1 \not = 0$. Let $u = {y_2 \over y_1}$.  Then $u \not = 0$.   The
derivation (\ref{e0}) still holds with $\uu = (1, u, u^2, \cdots,
u^{n-1})^\top \in \Re^n$.  Thus, $\PP$ is positive semi-definite.
\ep

\medskip

We may use the methods in \cite{QWW, QYW, QYX} to check if $\PP$ is
positive semi-definite or not when $m$ is even.   In Section 7, we
will present an algorithm for checking if $\PP$ is copositive or
not.

\medskip
\begin{Question} \label{Q1}
Can we give an example that $\PP$ is copositive but $\A$ is not?
\end{Question}

\begin{Question} \label{Q2}
When $m$ is even, can we give an example that $\PP$ is positive
semi-definite but $\A$ is not?
\end{Question}

\begin{Question} \label{Q3}
When $m$ is even, can we give an example that $\PP$ is positive
semi-definite but $\A$ is not?
\end{Question}

\begin{Question} \label{Q4}
Which conditions on $\PP$ may assure co-positiveness or positive
semi-definiteness of $\A$?
\end{Question}

\section{Strong Hankel Tensors and Generating Functions}

We are going to show that an even order strong Hankel tensor is
positive semi-definite.   In order to do this, we introduce a
generating function for a Hankel tensor $\A$.

Let $\A$ be a Hankel tensor defined by (\ref{e1}).    Let $f(t)$ be
an absolutely integrable real valued function on the real line
$(-\infty, \infty)$ such that
\begin{equation}\label{e2}
v_k \equiv \int_{-\infty}^\infty t^kf(t)dt, \end{equation} for $k =
0, \cdots, (n-1)m$.   Then we say that $f$ is a {\bf generating
function} of the Hankel tensor $\A$.    We see that $f(t)$ is also
the generating function of the associated Hankel matrix of $\A$. By
the theory of Hankel matrices \cite{Wik}, $f(t)$ is well-defined.

\begin{Theorem} \label{t1}
A Hankel tensor $\A$ has a nonnegative generating function if and
only if it is a strong Hankel tensor.  An even order strong Hankel
tensor is positive semi-definite.

On the other hand, suppose that $\A \in H_{m, n}$ has a generating
function $f(t)$ such that (\ref{e2}) holds. If $\A$ is copositive,
then
$$\int_{-\infty}^\infty t^{(i-1)m}f(t)dt \ge 0$$
for $i \in [n]$.
\end{Theorem}
\noindent {\bf Proof.}   By the famous Hamburger moment problem
\cite{Wik}, such a nonnegative generating function exists if and
only if the associated Hankel matrix is positive semi-definite,
i.e., $\A$ is a strong Hankel tensor.   On the other hand, suppose
that $\A$ has such a nonnegative generating function $f$ and $m$ is
even. Then for any $\x \in \Re^n$, we have
\begin{eqnarray*}
\A \x^{m} & = & \sum_{i_1, \cdots, i_m= 1}^n a_{i_1\cdots
i_m}x_{i_1}\cdots x_{i_m} \\
& = &\sum_{i_1, \cdots, i_m = 1}^n \int_{-\infty}^\infty
t^{i_1+\cdots +i_m-m}x_{i_1}\cdots x_{i_m}f(t)dt \\
& = & \int_{-\infty}^\infty \left(\sum_{i=1}^n
x_it^{i-1}\right)^mf(t)dt \\
& \ge  & 0. \end{eqnarray*} Thus, if $m$ is even and $\A$ is a
strong Hankel tensor, then $\A$ is positive semi-definite.

The final conclusion follows from (\ref{e2}) and Proposition
\ref{p1}. \ep

We now give an example of a positive semi-definite Hankel tensor,
which is not a strong Hankel tensor. Let $m=4$ and $n=2$. Let $v_0 =
v_4 = 1$, $v_2 = -{1 \over 6}$, and $v_1 = v_3 = 0$. Let $\A$ be
defined by (\ref{e1}).  Then for any $\x \in \Re^2$, we have
$$\A \x^4 = v_0 x_1^4 + 4 v_1 x_1^3x_2 + 6 v_2 x_1^2x_2^2 + 4
v_3x_1x_2^3 + v_4 x_2^4 = x_1^4 - x_1^2x_2^2 + x_2^4 \ge 0.$$ Thus,
$\A$ is positive semi-definite.   Let $A$ be the unique Hankel
matrix associated with $\A$.  Since $v_2 < 0$, by Proposition
\ref{p1}, $A$ is not positive semi-definite.    Thus, $\A$ is not a
strong Hankel tensor.

\medskip

\begin{Question} \label{Q5}
The question is, for a fixed even number $m \ge 4$, can we
characterize a positive semi-definite Hankel tensor by its
generating functions?
\end{Question}

\begin{Question} \label{Q6}
  If the associated Hankel matrix is
copositive, is the Hankel tensor copositive?
\end{Question}

\medskip

We now discuss the Hadamard product of two strong Hankel tensors.
Let $\A = (a_{i_1\cdots i_m}), \B = (b_{i_1\cdots i_m}) \in T_{m,
n}$. Define the Hadamard product of $\A$ and $\B$ as $\A \circ \B =
(a_{i_1\cdots i_m}b_{i_1\cdots i_m}) \in T_{m, n}$.  Clearly, the
Hadamard product of two Hankel tensors is a Hankel tensor.

\begin{Proposition} \label{p2}
The Hadamard product of two strong Hankel tensors is a strong Hankel
tensor.
\end{Proposition}
\noindent {\bf Proof.}   Let $\A$ and $\B$ be two strong Hankel
tensors in $H_{m, n}$.    Let $A$ and $B$ be Hankel matrices
associated with $\A$ and $\B$ respectively, such that $A$ and $B$
are positive semi-definite.    Clearly, the Hadamard product of $A$
and $B$ is a Hankel matrix associated with the Hadamard product of
$\A$ and $\B$. By the Schur product theorem \cite{HJ}, the Hadamard
product of two positive semi-definite symmetric matrices is still a
positive semi-definite symmetric matrix.   Thus, the Hadamard
product of $A$ and $B$ is positive semi-definite.  This implies that
the Hadamard product of $\A$ and $\B$ is a strong Hankel tensor. \ep

On the other hand, the Hadamard product of two positive
semi-definite Hankel tensors may not be positive semi-definite.
Assume that $m=4$ and $n=2$. Let $\A$ be the example given above.
Then $\A$ is a positive semi-definite Hankel tensor. On the other
hand, let $\B = (b_{i_1i_2i_3i_4}) \in S_{4, 2}$ be defined by
$b_{i_1i_2i_3i_4} = 1$ if $i_1+i_2+i_3+i_4 = 6$, and
$b_{i_1i_2i_3i_4} = 0$ otherwise.  We may verify that $\B$ is a
strong Hankel tensor, thus a positive semi-definite Hankel tensor.
It is easy to verify that $\A \circ \B$ is not positive
semi-definite.    Note here that $\A$ is not a strong Hankel tensor.
Thus, this example does not contradict Proposition \ref{p2}.

\section{Vandermonde Decomposition and Complete Hankel Tensors}

For any vector $\uu \in \Re^n$, $\uu^m$ is a rank-one $m$th order
symmetric $n$-dimensional tensor $\uu^m = \left(u_{i_1}\cdots
u_{i_m}\right) \in S_{m, n}$.   If $\uu = (1, u, u^2, \cdots,
u^{n-1})^\top$, then $\uu$ is called a {\bf Vandermonde vector}
\cite{PDV}. If
\begin{equation}\label{e3}
\A = \sum_{k=1}^r \alpha_k \left(\uu_k\right)^m,
\end{equation}
where $\alpha_k \in \Re$, $\alpha_k \not = 0$, $\uu_k = (1, u_k,
u_k^2, \cdots, u_k^{n-1})^\top\in \Re^n$ are Vandermonde vectors for
$k= 1, \cdots, r$, and $u_i \not = u_j$ for $i \not = j$, then we
say that tensor $\A$ has a {\bf Vandermonde decomposition}. We call
the minimum value of $r$ the {\bf Vandermonde rank} of $\A$.

\begin{Theorem} \label{t2}
Let $\A \in S_{m,n}$.  Then $\A$ is a Hankel tensor if and only if
it has a Vandermonde decomposition (\ref{e3}).   In this case, we
have $r \le (n-1)m+1$.

Suppose that $\A$ has a Vandermonde decomposition (\ref{e3}).   If
$\A$ is copositive, then
\begin{equation}\label{e3.5}
\sum_{k=1}^r \alpha_k u_k^{(i-1)m} \ge 0, \ \ {\rm for}\ i \in [n].
\end{equation}
On the other hand, if $m$ is even and $\alpha_k > 0$ for $i \in
[r]$, then $\A$ is positive semi-definite.
\end{Theorem}

\noindent {\bf Proof.}   Suppose that $\A$ has a Vandermonde
decomposition (\ref{e3}).   Let
\begin{equation}\label{e4}
v_i = \sum_{k=1}^r \alpha_k u_k^i, \ \ {\rm for}\ i = 0, \cdots,
(n-1)m.
\end{equation}
By (\ref{e3}), we see that (\ref{e1}) holds.   Thus, $\A$ is a
Hankel tensor.

On the other hand, assume that $\A$ is a Hankel tensor defined by
(\ref{e1}).    Let $r = (n-1)m+1$. Pick real numbers $u_k, k \in
[r]$ such that $u_i \not = u_j$ for $i \not = j$. By matrix analysis
\cite{HJ}, the coefficient matrix of the linear system (\ref{e4})
with $\alpha_k, k \in [r]$ as variables, is a Vandermonde matrix,
which is nonsingular.   Thus, the linear system (\ref{e4}) has a
solution $\alpha_k, k \in [r]$. Substituting such $\alpha_k, k = 1,
\cdots, r$ to (\ref{e3}), we see that (\ref{e3}) holds, i.e., $\A$
has a Vandermonde decomposition.

Suppose that $\A$ has a Vandermonde decomposition (\ref{e3}).   If
$\A$ is copositive, then (\ref{e3.5}) follows from (\ref{e4}) and
Proposition \ref{p1}.   On the other hand, assume that $m$ is even.
Suppose (\ref{e3}) holds with $\alpha_k
> 0, k \in [r]$.  For any $\x \in \Re^n$, we have
$$\A \x^m = \sum_{k=1}^r \alpha_k (\uu_k^\top \x)^m \ge 0.$$
Thus, $\A$ is positive semi-definite. \ep

\medskip

In (\ref{e3}), if $\alpha_k > 0, k \in [r]$, then we say that $\A$
has a positive Vandermonde decomposition and call $\A$ a {\bf
complete Hankel Tensor}.   Thus, Theorem \ref{t2} says that an even
order complete Hankel tensor is positive semi-definite.   We will
study the spectral properties of odd order complete Hankel tensors
in the next section.

By (\ref{e4}), if $\alpha_k > 0$ for $k \in [r]$, then $v_i$ is
nonnegative if $i$ is even.    Thus, the counterexample $\A$, given
in the last section, is not a complete Hankel tensor as it has $v_2
< 0$.  This implies that a positive semi-definite Hankel tensor may
not be a complete Hankel tensor.

\medskip
We now discuss the Hadamard product of two complete Hankel tensors.

\begin{Proposition} \label{p2.5}
The Hadamard product of two complete Hankel tensors is a complete
Hankel tensor.
\end{Proposition}
\noindent {\bf Proof.}  Suppose that $\A, \B \in H_{m, n}$ are two
complete Hankel tensors.   Then we may assume that each of $\A$ and
$\B$ has a positive Vandermonde decomposition:
$$\A = \sum_{k=1}^r \alpha_k \left(\uu_k\right)^m$$
and
$$\B = \sum_{j=1}^s \beta_j \left(\vv_j\right)^m,$$
where $\alpha_k > 0, \uu_k = (1, u_k, u_k^2, \cdots,
u_k^{n-1})^\top$ are Vandermonde vectors for $k \in [r]$, $\beta_j
> 0, \vv_j = (1, v_j, v_j^2, \cdots, v_j^{n-1})^\top$ are
Vandermonde vectors for $j \in [s]$. Then the Vandermonde product of
$\A$ and $\B$ is
$$\A \circ \B = \sum_{k=1}^r \sum_{j=1}^s \alpha_k\beta_j
\left(\w_{kj}\right)^\top,$$ where $\alpha_k\beta_j > 0, \w_{kj}
=(1, u_kv_j, (u_kv_j)^2, \cdots, (u_kv_j)^{n-1})^\top$ are
Vandermonde vectors for $k \in [r]$ and $j \in [s]$.   We see that
$\A \circ \B$ has a positive Vandermonde decomposition, thus a
complete Hankel tensor.   \ep

We may summarize the results on Hadamard products.   The Hadarmard
product of two Hankel tensors is a Hankel tensor.  The Hadarmard
product of two strong Hankel tensors is a strong Hankel tensor. The
Hadarmard product of two complete Hankel tensors is a complete
Hankel tensor. But the Hadarmard product of two positive
semi-definite Hankel tensors may not be positive semi-definite.

\medskip

\begin{Question} \label{Q7}
Can we characterize a positive semi-definite Hankel tensor by its
Vandermonde decomposition?
\end{Question}

\begin{Question} \label{Q8}
Is a strong Hankel tensor a complete Hankel tensor?  Is a complete
Hankel tensor a strong Hankel tensor?
\end{Question}

\section{Spectral Properties of Odd Order Complete and Strong Hankel Tensors}

Suppose that $m$ is even.  Then by Theorem 5 of \cite{Qi}, all the
H-eigenvalues and Z-eigenvalues of a strong Hankel tensor or a
complete Hankel tensor are nonnegative, as strong Hankel tensors and
complete Hankel tensors are positive semi-definite.  In this
section, we discuss spectral properties of odd order complete and
strong Hankel tensors.   Hence, assume that $m$ is odd in this
section.

We now briefly review the definition of eigenvalues, H-eigenvalues
E-eigenvalues and Z-eigenvalues of a real $m$th order
$n$-dimensional symmetric tensor $\A = \left(a_{i_1\cdots
i_m}\right) \in S_{m, n}$ \cite{Qi}.    Let $\x = (x_1, \cdots,
x_n)^\top \in {\boldmath C}^n$. Then $\A \x^{m-1}$ is an
$n$-dimensional vector, with its $i$th component as \noindent
$\sum_{i_2\cdots i_m=1}^n a_{ii_2\cdots i_m}x_{i_2}\cdots x_{i_m}$.
For any vector $\x \in {\boldmath C}^n$, $\x^{[m-1]}$ is a vector in
${\boldmath C}^n$, with its $i$th component as $x_i^{m-1}$.  If $\A
\x^{m-1} = \lambda \x^{[m-1]}$ for some $\lambda \in {\boldmath C}$
and $\x \in {\boldmath C}^n \setminus \{ 0 \}$, then $\lambda$ is
called an {\bf eigenvalue} of $\A$ and $\x$ is called an {\bf
eigenvector} of $\A$, associated with $\lambda$. If both $\lambda$
and $\x$ are real, then they are called an {\bf H-eigenvalue} and an
{\bf H-eigenvector} of $\A$, respectively. If $\A x^{m-1} = \lambda
\x$ for some $\lambda \in {\boldmath C}$ and $\x \in {\boldmath
C}^n$, satisfying $\x^\top \x = 1$, then $\lambda$ is called an {\bf
E-eigenvalue} of $\A$ and $\x$ is called an {\bf E-eigenvector} of
$\A$, associated with $\lambda$. If both $\lambda$ and $\x$ are
real, then they are called a {\bf Z-eigenvalue} and a {\bf
Z-eigenvector} of $\A$, respectively. Note \cite{Qi} that
Z-eigenvalues always exist, and when $m$ is even, H-eigenvalues
always exist.

\begin{Proposition} \label{p3}
Suppose that $m$ is odd and $\A \in H_{m, n}$ is a complete Hankel
tensor.  Assume that $\A$ has at least one H-eigenvalue.   Then all
the H-eigenvalues of $\A$ are nonnegative. Let $\lambda$ be an
H-eigenvalue of $\A$, with an H-eigenvector $\x = (x_1, \cdots,
x_n)^\top$.   Then either $\lambda = 0$ or $\lambda > 0$ with $x_1
\not = 0$.
\end{Proposition}
\noindent {\bf Proof.}  By the definition of complete Hankel
tensors, $\A$ has a Vandermonde decomposition (\ref{e3}), with
$\alpha_k > 0$ for $k \in [r]$.  Suppose that $\A$ has an
H-eigenvalue $\lambda$ associated with an H-eigenvector $\x = (x_1,
\cdots, x_n)^\top$.   Then for $i \in [n]$, we have
\begin{equation} \label{e5}
\lambda x_i^{m-1} = \left(\A \x^{m-1}\right)_i = \sum_{k=1}^r
\alpha_k u_k^{i-1} \left[(\uu_k)^\top \x\right]^{m-1}.
\end{equation}
If $(\uu_k)^\top \x = 0$ for all $k \in [r]$, then the right hand
side of (\ref{e5}) is $0$.   Since $\x \not = \0$, we may pick $i$
such that $x_i \not = 0$. Then (\ref{e5}) implies that $\lambda =
0$.

Suppose that $(\uu_k)^\top \x \not = 0$ for at least one $k$.   Let
$i=1$.   Then the the right hand side of (\ref{e5}) is positive.
This implies that $\lambda > 0$ and $x_1 \not = 0$.   \ep

In general an odd order symmetric tensor may not have H-eigenvalues.

\begin{Question} \label{Q9}
Does a complete Hankel tensor always have an H-eigenvalue?
\end{Question}

\medskip

For Z-eigenvalues, we have the following results.

\begin{Proposition} \label{p4}
Suppose that $m$ is odd and $\x = (x_1, \cdots, x_n)^\top$ is a
Z-eigenvector of a complete Hankel tensor $\A \in H_{m, n}$,
associated with a Z-eigenvalue $\lambda$.    Then $x_i \ge 0$ for
all odd $i$ and $x_1
> 0$ if $\lambda
> 0$; and $x_i \le 0$ for all odd $i$ and $x_1 < 0$ if $\lambda < 0$.
\end{Proposition}
\noindent {\bf Proof.}  Again, by the definition of complete Hankel
tensors, $\A$ has a Vandermonde decomposition (\ref{e3}), with
$\alpha_k > 0$ for $k \in [r]$.  Suppose that $\A$ has a
Z-eigenvalue $\lambda$ associated with a Z-eigenvector $\x = (x_1,
\cdots, x_n)^\top$. Then for $i \in [n]$, we have
\begin{equation} \label{e6}
\lambda x_i = \left(\A \x^{m-1}\right)_i = \sum_{k=1}^r \alpha_k
u_k^{i-1} \left[(\uu_k)^\top \x\right]^{m-1}.
\end{equation}
If $(\uu_k)^\top \x = 0$ for all $k \in [r]$, then the right hand
side of (\ref{e6}) is $0$.   Since $\x \not = \0$, we may pick $i$
such that $x_i \not = 0$.   Then (\ref{e6}) implies that $\lambda =
0$.

Suppose that $(\uu_k)^\top \x \not = 0$ for at least one $k$.   Let
$i$ be odd.   Then the the right hand side of (\ref{e6}) is
nonnegative. This implies that $\lambda x_i \ge 0$.  The conclusion
on $x_i$ with $i$ odd follows.  Let $i=1$.   Then the the right hand
side of (\ref{e6}) is positive. This implies that $\lambda x_1
> 0$.  The conclusion on $x_1$ follows now. \ep

\medskip

We now study spectral properties of odd order strong Hankel tensors.

\begin{Proposition} \label{p5}
Suppose that $m$ is odd and $\x = (x_1, \cdots, x_n)^\top$ is a
Z-eigenvector of a strong Hankel tensor $\A \in H_{m, n}$,
associated with a Z-eigenvalue $\lambda$.    Then $x_i \ge 0$ for
all odd $i$ if $\lambda
> 0$; and $x_i \le 0$ for all odd $i$ if $\lambda < 0$.
\end{Proposition}
\noindent {\bf Proof.}  By Theorem \ref{t1}, $\A$ has a nonnegative
generating function $f(t)$ such that (\ref{e2}) holds. Suppose that
$\A$ has a Z-eigenvalue $\lambda$ associated with a Z-eigenvector
$\x = (x_1, \cdots, x_n)^\top$. Then for $i \in [n]$, we have
\begin{eqnarray} \label{e7}
\lambda x_i & = & \left(\A \x^{m-1}\right)_i  \nonumber \\
& =  & \sum_{i_2, \cdots, i_m= 1}^n a_{ii_2\cdots i_m}x_{i_2}\cdots
x_{i_m} \nonumber \\
& = & \sum_{i_2, \cdots, i_m = 1}^n \int_{-\infty}^\infty t^{i+
i_2+\cdots +i_m-m}x_{i_2}\cdots x_{i_m}f(t)dt \nonumber \\
& = & \int_{-\infty}^\infty t^{i-1}\left(\sum_{j=1}^n
x_jt^{j-1}\right)^{m-1}f(t)dt.
\end{eqnarray}
Let $i$ be odd.   Then the the right hand side of (\ref{e7}) is
nonnegative.  The conclusion follows now. \ep

Note that we miss a result of the H-eigenvalues of an odd order
strong Hankel tensor.

\begin{Question} \label{Q10}
Are all the H-eigenvalues of an odd order strong Hankel tensor
nonnegative?
\end{Question}

Similar spectral properties hold for odd order Laplacian tensors
\cite{Qi1} and odd order completely positive tensors \cite{QXX}. A
common point is that such classes of symmetric tensors are positive
semi-definite when the order is even.   Thus, we may think if we may
define some odd order ``positive semi-definite'' symmetric tensors,
with such spectral properties.   Further study is needed on such a
phenomenon.

\section{Upper Bounds for the Smallest Z-Eigenvalue and Lower
Bounds for the Largest Z-Eigenvalue}

Let $\A \in S_{m, n}$.  Then $\A$ always has Z-eigenvalues
\cite{Qi}.   Denote the smallest and the largest Z-eigenvalue of
$\A$ by $\lambda_{\min}(\A)$ and $\lambda_{\max}(\A)$ respectively.
We always have \cite{Qi}
\begin{equation} \label{e10}
\lambda_{\min}(\A) = \min \{ \A \x^m : \x \in \Re^n, \x^\top \x = 1
\}
\end{equation}
and
\begin{equation} \label{e11}
\lambda_{\max}(\A) = \max \{ \A \x^m : \x \in \Re^n, \x^\top \x = 1
\}.
\end{equation}

 If $m$ is even, $\A$ is positive semi-definite if and only if
$\lambda_{\min}(\A) \ge 0$ \cite{Qi}.   If $m$ is odd, then
$\lambda_{\max}(\A) \ge 0$ and $\lambda_{\min}(\A) = -
\lambda_{\max}(\A)$.   In general, $\max \{ |\lambda_{\min}(\A)|,
|\lambda_{\max}(\A)| \}$ is a norm of $\A$ in the space $S_{m, n}$
\cite{QiL}.  If $|\lambda_{\min}(\A)| = \max \{
|\lambda_{\min}(\A)|, |\lambda_{\max}(\A)| \}$, then
$\lambda_{\min}(\A)$ and its corresponding eigenvector $\x$ form the
best rank-one approximation to $\A$ \cite{Qi, QWW}.  Similarly, if
$|\lambda_{\max}(\A)| = \max \{ |\lambda_{\min}(\A)|,
|\lambda_{\max}(\A)| \}$, then $\lambda_{\max}(\A)$ and its
corresponding eigenvector $\x$ form the best rank-one approximation
to $\A$ \cite{Qi, QWW}.   Let $\x \in \Re^n, \x \not = \0$.   By
(\ref{e10}) and (\ref{e11}), we have
\begin{equation} \label{e12}
\lambda_{\min}(\A)  \le {\A \x^m \over \| \x \|_2^m} \le
\lambda_{\max}(\A).
\end{equation}

With the above knowledge, for a Hankel tensor $\A$, we may give some
upper bounds for $\lambda_{\min}(\A)$,  and some lower bounds for
$\lambda_{\max}(\A)$.

\begin{Proposition} \label{p6}
Suppose that $\A \in H_{m, n}$.  Then
$$\lambda_{\min}(\A) \le \min_{i \in [n]} v_{(i-1)m} \le \max_{i \in [n]}
v_{(i-1)m} \le \lambda_{\max}(\A).$$
\end{Proposition}
\noindent {\bf Proof.} Since $v_{(i-1)m} = \A(\e_i)^m$ for $i \in
[n]$, the conclusion follows from (\ref{e12}). \ep.

Suppose $\PP$ is the associated plane tensor of $\A$.  We now use
$\lambda_{\min}(\PP)$ and $\lambda_{\max}(\PP)$ to give an upper
bound for $\lambda_{\min}(\A)$, and a lower bound for
$\lambda_{\max}(\A)$, respectively.

\begin{Proposition} \label{p7}
Suppose that $\A \in H_{m, n}$, and $\PP$ is the associated plane
tensor of $\A$.  Assume that $m(n-1)$ is even.  If $\y = (y_1,
y_2)^\top$ is a Z-eigenvector of $\PP$, associated with
$\lambda_{\min}(\PP)$
\begin{equation} \label{e13}
\sqrt{\sum_{j=0}^{(n-1)m} y_1^{2(n-1)m-2j}y_2^{2j}}
\lambda_{\min}(\A) \le \lambda_{\min}(\PP).
\end{equation}
If $\z = (z_1, z_2)^\top$ is a Z-eigenvector of $\PP$, associated
with $\lambda_{\max}(\PP)$
\begin{equation} \label{e14}
\sqrt{\sum_{j=0}^{(n-1)m} z_1^{2(n-1)m-2j}z_2^{2j}}
\lambda_{\max}(\A)    \ge \lambda_{\max}(\PP).
\end{equation}
\end{Proposition}
\noindent {\bf Proof.}  If $y_1 = 0$, since $y_1^2 + y_2^2 = 1$,
then
$$\sqrt{\sum_{j=0}^{(n-1)m}
y_1^{2(n-1)m-2j}y_2^{2j}} = 1.$$ We have $$ \lambda_{\min}(\PP) =
\PP \y^{(n-1)m} = v_{(n-1)m} \ge \lambda_{\min}(\A) =
\sqrt{\sum_{j=0}^{(n-1)m} y_1^{2(n-1)m-2j}y_2^{2j}}
\lambda_{\min}(\A),$$ where the inequality is due to Proposition
\ref{p6}.   Thus, (\ref{e13}) holds.

Suppose that $y_1 \not = 0$. Let $u = {y_2 \over y_1}$ and $\uu =
(1, u, u^2, \cdots, u^{n-1})^\top \in \Re^n$. Then
\begin{eqnarray*}
\lambda_{\min}(\PP) & =  & \PP \y^{(n-1)m} \\
 & = & y_1^{(n-1)m} \sum_{k=0}^{(n-1)m} \left({(n-1)m
\over k}\right) \cdot { s_{k, m} v_k \over \left({(n-1)m \atop
k}\right)} u^k \\
& = & \left|y_1^{(n-1)m}\right| \A \uu^m \\
& = & \left|y_1^{(n-1)m}\right| \|\uu\|_2^m {\A \uu^m \over \|\uu\|_2^m } \\
& = & \sqrt{\sum_{j=0}^{(n-1)m} y_1^{2(n-1)m-2j}y_2^{2j}} {\A \uu^m
\over \|\uu\|_2^m} \\
& \ge & \sqrt{\sum_{j=0}^{(n-1)m} y_1^{2(n-1)m-2j}y_2^{2j}}
\lambda_{\min}(\A),
\end{eqnarray*}
where the inequality is due to (\ref{e12}).   Thus, (\ref{e13}) also
holds in this case.   This proves (\ref{e13}).

We may prove (\ref{e14}) similarly. \ep

\begin{Question} \label{Q11}
Suppose that a Hankel tensor $\A$ is associated with a Hankel matrix
$A$.   Can we use the largest and the smallest eigenvalues of $A$ to
bound the largest and the smallest H-eigenvalues (Z-eigenvalues) of
$\A$?
\end{Question}

\section{An Algorithm for Recognizing Copositivity of  a Symmetric Plane Tensor}

In Section 2, we show that if a Hankel tensor $\A \in H_{m, n}$ is
copositive, then its associated plane tensor $\PP \in S_{(n-1)m, 2}$
must be copositive.    In this section, we present an algorithm to
determine a plane tensor $\PP \in S_{l, 2}$ is copositive or not.
Here, $l \ge 2$.

Let $\PP = (p_{i_1\cdots i_l})$.   Denote $p_k = p_{i_1\cdots i_l}$
if $k$ of $i_1, \cdots, i_l$ are $2$ and the others are $1$.   Then
for any $\y = (y_1, y_2)^\top \in \Re^2$, we have
$$\PP \y^l = \sum_{k=0}^l \left({l \atop k} \right) p_ky_1^{l-k}y_2^k.$$
It is easy to see that $\PP$ is copositive if and only if
$$\min \{ \PP \y^l : y_1 + y_2 = 1, y_1 \ge 0, y_2 \ge 0 \} \ge 0,$$
i.e.,
\begin{equation} \label{e15}
\min \left\{  \sum_{k=0}^l \left({l \atop k} \right)
p_ky_1^{l-k}y_2^k : y_1 + y_2 = 1, y_1 \ge 0, y_2 \ge 0 \right\} \ge
0.
\end{equation}
Let $t = y_1$.   Then $y_2 = 1-t$.   We may rewrite (\ref{e15}) as
\begin{equation} \label{e16}
\min \left\{ \phi(t) : 0 \le t \le 1 \right\} \ge 0,
\end{equation}
where
\begin{equation} \label{e17}
\phi(t) = \sum_{k=0}^l \left({l \atop k} \right) p_kt^{l-k}(1-t)^k.
\end{equation}
To check if (\ref{e16}) holds, we only need to check if $\phi(t) \ge
0$ for all critical points $t$ of (\ref{e16}).  By optimization
theory, the critical points of (\ref{e16}) are $t=0$, $t=1$, and any
$t \in (0, 1)$ such that $\phi'(t) = 0$.   Note that $\phi(0) = p_l$
and $\phi(1) = p_0$.    Thus, we have a simple algorithm to check if
$\PP$ is copositive or not.

\begin{Algorithm}

\noindent
\smallskip

Step 1. If $p_0 < 0$ or $p_l < 0$, then $\PP$ is not copositive.
Stop. Otherwise, go to the next step.

\smallskip

Step 2.  Find all the critical points $t$ such that $\phi'(t) = 0$
and $0 < t < 1$, where $\phi(t)$ is defined by (\ref{e17}).  If
$\phi(t) < 0$ for one of such critical point $t$, then $\PP$ is not
copositive.  Otherwise $\PP$ is copositive.  Stop.
\end{Algorithm}

We see that this algorithm is simple.

\section{Final Remarks and Further Questions}

In this paper, we make an initial study on Hankel tensors.   We see
that Hankel tensors have a very special structure, hence have very
special properties.  We associate a Hankel tensor with a Hankel
matrix, a symmetric plane tensor, generating functions and
Vandermonde decompositions.   They will be useful tools for further
study on Hankel tensors.

Some questions have already been raised in Sections 2-6. Here are
some further questions.

\begin{Question} \label{Q12}
Badeau and Boyer \cite{BB} proposed fast higher-order singular value
decomposition (HOSVD) for third order Hankel tensors.   Can we
construct some efficient algorithms for the largest and the smallest
H-eigenvalues (Z-eigenvalues) of a Hankel tensor, or a strong Hankel
tensor, or a complete Hankel tensor?
\end{Question}

\begin{Question} \label{Q13}
In general, it is NP-hard to compute the largest and the smallest
H-eigenvalues (Z-eigenvalues) of a symmetric tensor.  What is the
complexity for computing the smallest H-eigenvalues (Z-eigenvalues)
of a Hankel tensor, a strong Hankel tensor, and a complete Hankel
tensor?
\end{Question}

\begin{Question} \label{Q14}
Proposition 8 of \cite{Qi} says that the determinants of all the
principal symmetric sub-tensors of a positive semi-definite tensor
are nonnegative.    The converse is not true in general.   Is the
converse of Proposition 8 of \cite{Qi} true for Hankel tensors?
\end{Question}

For the definition of the determinants of tensors, see \cite{HHLQ,
Qi, SSZ}. They were called symmetric hyperdeterminants in \cite{Qi},
and simply determinants in \cite{HHLQ, SSZ}.

\begin{Question} \label{Q15}
The theory of Hankel matrices is based upon finite and infinite
Hankel matrices as well as Hankel operators \cite{Wid}.   Should we
also study infinite Hankel tensors and multi-linear Hankel
operators?
\end{Question}

\bigskip

\begin{center}
{\bf Acknowledgment}
\end{center}

The author is thankful to Professor Changqing Xu, who suggested the
research topic ``Hankel tensors'', to Mr. Zhongming Chen, who
suggested the proof of Proposition \ref{p2.5}, to Dr. Yisheng Song,
who made some comments, and to two referees for their comments.


\end{document}